\documentclass[12pt]{article}
\usepackage{amsmath}
\usepackage{amssymb}
\usepackage{a4wide}
\newtheorem{thm}{Theorem}[section]
\newtheorem{la}[thm]{Lemma}
\newtheorem{Defn}[thm]{Definition}
\newtheorem{Exam}[thm]{Example}
\newtheorem{Remark}[thm]{Remark}

\newenvironment{rem}{\begin{Remark}\rm}{\end{Remark}}
\newenvironment{proof}{{\noindent\bf Proof.}}%
                  {\nopagebreak\hspace*{\fill}$\Box$\medskip\medskip\par}   
\newcommand{\Punkt}{\nopagebreak\hspace*{\fill}$\Box$}

\newcommand{\mto}{\mapsto}

\newcommand{\N}{{\mathbb N}}
\newcommand{\Z}{{\mathbb Z}}
\newcommand{\R}{{\mathbb R}}
\newcommand{\cS}{{\mathcal S}}
\newcommand{\cK}{{\mathcal K}}
\newcommand{\dl}{{\displaystyle \lim_{\longrightarrow}}}
\newcommand{\sub}{\subseteq}

\newcommand{\wb}{\overline}
\DeclareMathOperator{\Ad}{Ad}
\DeclareMathOperator{\Hol}{Hol}
\DeclareMathOperator{\Exp}{Exp}
\begin{document}
\begin{center}
{\Large\bf Direct limit groups do not have small
subgroups}\vspace{4.5 mm}\\
{\bf Helge Gl\"{o}ckner}
\end{center}
\begin{abstract}
\noindent
We show that countable direct limits of
finite-dimensional Lie groups do not have small
subgroups. The same conclusion is
obtained for suitable direct limits
of infinite-dimensional Lie groups.\vspace{-1.1mm}
\end{abstract}
\begin{center}
{\bf\Large Introduction}\vspace{-.4mm}
\end{center}
The present investigation is related to an open problem in the
theory of infinite-dimensional Lie groups, i.e., Lie groups
modelled on locally convex spaces (as in~\cite{Mil}).
Recall that a topological group~$G$ is said to \emph{have small
subgroups} if every identity neighbourhood
$U\sub G$ contains a non-trivial subgroup of~$G$.
If every identity neighbourhood $U$ contains a non-trivial
torsion group, then $G$ is said to \emph{have small torsion
subgroups}.
The additive group of the Fr\'{e}chet space $\R^\N$ is an example
of a Lie group which has small subgroups.
It is an open problem (formulated first in \cite{Nee}) whether
a Lie group modelled on a locally convex space can have
small torsion subgroups.
As a general proof for the non-existence of small torsion
subgroups seems to be out of reach, it is natural to examine
at least the main examples of infinite-dimensional Lie groups,
and to rule out this pathology individually for each of them.
The main examples comprise linear Lie groups, diffeomorphism groups,
mapping groups, and \emph{direct limit groups},
i.e., direct limits in the
category of Lie groups of countable direct systems of
finite-dimensional Lie groups, as constructed in \cite{FUN}
(see also \cite{DIR}, \cite[Theorem~47.9]{KaM} and \cite{NRW}
for special cases).
We show that direct limit
groups do not have small subgroups,
thus ruling out the existence of
small torsion subgroups in particular:\\[3.1mm]
{\bf Theorem~A.}
{\em Let $\cS:=((G_n)_{n\in \N}, (i_{n,m})_{n\geq m})$
be a direct sequence of finite-dimensional
real Lie groups $G_n$ and smooth homomorphisms
$i_{n,m}\colon G_m\to G_n$.
Let $G=\dl\,G_n$\vspace{-.7mm}
be the direct limit of~$\cS$
in the category of Lie groups
modelled on locally convex spaces.
Then $G$ does not have small
subgroups.}\\[3.1mm]
More generally, we can tackle
direct limits of
not necessarily finite-dimensional
Lie groups.\\[3.1mm]
{\bf Theorem~B.}
{\em Let $G$ be a Lie group
modelled on a locally convex space
which is the union
of an ascending sequence $G_1\leq G_2\leq\cdots$
of Lie groups~$G_n$ modelled
on locally convex spaces,
such that the inclusion maps
$i_{n,m}\colon G_m\to G_n$
for $m\leq n$
and $i_n\colon G_n\to G$
are smooth homomorphisms.
Assume that at least one of the following
conditions is satisfied:
\begin{itemize}
\item[\rm (i)]
Each $G_n$ is a Banach-Lie group,
$L(i_{n,m})\colon L(G_m)\to L(G_n)$
is a compact operator for all positive integers
$m<n$, and $G=\dl\,G_n$\vspace{-2.4mm}
as a topological space; or:
\item[\rm (ii)]
$G$ admits a direct limit chart,
$L(G_n)$ is a $k_\omega$-space
admitting a continuous norm,
and $G_n$ has an exponential map
which is a local homeomorphism at~$0$,
for each $n\in \N$.
\end{itemize}
Then $G$ does not
have small subgroups.}\pagebreak

\noindent
{\bf Remarks.}
\begin{itemize}
\item[(a)]
All of the maps $L(i_{n,m})$
are injective in Theorem~B,
since $i_{n,m}$ is an injective smooth homomorphism
and $G_m$ has an exponential function
(cf.\ \cite[Lemma~7.1]{Mil}).
\item[(b)]
A Hausdorff topological space~$X$ is called
a \emph{$k_\omega$-space} if there exists an ascending
sequence $K_1\sub K_2\sub \cdots$ of compact subsets
of~$X$ such that $X=\bigcup_{n\in \N}K_n$
and $U\sub X$ is open if and only if
$U\cap K_n$ is open in~$K_n$, for each $n\in \N$
(i.e., $X=\dl\,K_n$\vspace{-.6mm} as a topological space).
Then $(K_n)_{n\in \N}$ is called a
\emph{$k_\omega$-sequence} for~$X$.
For background information concerning
$k_\omega$-spaces with a view towards direct limit
constructions, see \cite{GaG} and the references therein.
\item[(c)]
A locally convex space~$E$
is a \emph{Silva space}
(or ($LS$)-space)
if it is the locally convex direct limit
$E=\bigcup_{n\in \N}E_n=\dl\,E_n$\vspace{-.3mm} of a
sequence $E_1\sub E_2\sub \cdots$
of Banach spaces and each inclusion map $E_n\to E_{n+1}$
is a compact linear operator.
Then $E=\dl\,E_n$\vspace{-.4mm}
as a topological space \cite[\S7.1, Satz]{Flo},
and $E$ is a $k_\omega$-space~\cite[Example~9.4]{COM}.
It is also known that the dual space
$E'$ of any
metrizable locally convex space~$E$
is a $k_\omega$-space, when equipped with the
topology of compact convergence
(cf.\ \cite[Corollary~4.7]{Aus}).
\item[(d)]
By definition,
the existence of a \emph{direct limit chart}
means the following:
$L(G)=\dl\, L(G_n)$\vspace{-.3mm}
as a locally convex space,
and there exists a chart
$\phi\colon U\to V\sub L(G)$ of~$G$ around~$1$,
with the following properties:
$U=\bigcup_{n\in \N}U_n$,
$V=\bigcup_{n\in \N}V_n$
and $\phi=\bigcup_{n\in \N}\phi_n=\dl\,\phi_n$\vspace{-.7mm}
for certain charts
$\phi_n\colon U_n\to V_n\sub L(G_n)$
of~$G_n$ around~$1$,
satisfying $U_n\sub U_{n+1}$
and $\phi_{n+1}|_{U_n}=\phi_n$
for each $n\in \N$
(see \cite{COM} for further
information).
\item[(e)]
For example, every direct limit
of an ascending sequence of
finite-dimensional Lie groups
admits a direct limit chart,
by construction of the Lie group
structure in~\cite{FUN}.
In the situation of Theorem~A,
we may always assume that each $i_{n,m}$
(and hence also each limit map
$i_n\colon G_n\to G$) is injective
(see \cite[Theorem~4.3]{FUN}).
Then $G=\dl\,G_n$\vspace{-.6mm}
as a topological space
by \cite[Theorem~4.3\,(a)]{FUN}.
Thus Theorem~A is a special case
of Theorem~B\,(i) and does not require a separate
proof.
\item[(f)] If condition\,(ii) of
Theorem~B is satisfied,
then $L(G)$ is a $k_\omega$-space
and $L(G)=\dl\, L(G_n)$\vspace{-.3mm}
as a topological space,
by \cite[Proposition~7.12]{GaG}.
If $\phi\colon U\to V\sub L(G)$
is a direct limit chart for~$G$,
with $U=\bigcup_{n\in \N}U_n$
and $V=\bigcup_{n\in \N}V_n$
as in\,(d),
then
$V=\dl\,V_n$\vspace{-.3mm}
as a topological space (see Lemma~\ref{baseDL}\,(b) below)
and hence $U=\dl\,U_n$.\vspace{-.3mm}
Using translations, it easily follows that
also $G=\dl\,G_n$\vspace{-.3mm}
as a topological space.
\item[(g)]
Suppose that $G_n$ is a Banach-Lie group
in the situation of Theorem~B,
$L(i_{n,m})$ is a compact operator
for $n>m$, and $G$
admits a direct limit chart.
Then $G=\dl\,G_n$\vspace{-.5mm}
as a topological space
(since (c) allows us to repeat the argument from\,(f)),
and thus condition\,(i) of Theorem~B
is satisfied. While the direct limit
property required in\,(i)
is somewhat elusive, the existence
of a direct limit chart can frequently
be checked in concrete situations.
\end{itemize}
{\bf Example.}
To illustrate the use of Theorem~B\,(i),
let $H$ be a finite-dimensional complex
Lie group and $K$ be a non-empty compact subset
of a finite-dimensional complex vector space~$X$.
Then the group $\Gamma(K,H)$
of germs of complex analytic
$H$-valued maps on open neighbourhoods
of~$K$ is a Lie group in a natural way.
It is modelled on the locally convex direct limit
$\Gamma(K,L(H))=\dl\, \Hol_b(U_n,L(H))$,\vspace{-.3mm}
where $U_1\supseteq U_2\supseteq\cdots$
is a fundamental sequence of open
neighbourhoods of~$K$
with $U_{n+1}$ relatively compact
in~$U_n$, for each $n\in \N$,
and such that each connected component
of $U_n$ meets~$K$.
Furthermore, $\Hol_b(U_n,L(H))$ denotes
the Banach space of bounded holomorphic
functions from $U_n$ to~$L(H)$,
equipped with the supremum norm.
For the identity component, we have
$G:=\Gamma(K,H)_0=\dl\,G_n$\vspace{-.6mm}
for certain Banach-Lie groups~$G_n$
satisfying condition\,(i)
of Theorem~B, and thus $G$ does not have
small subgroups (nor $\Gamma(K,H)$).\\[2.8mm]
In fact, let
$\Hol(U_n,H)$ be the group
of all complex analytic
$H$-valued maps on~$U_n$.
Since $\Exp_n\colon \Hol_b(U_n,L(H))\to \Hol(U_n,H)$,
$\Exp_n(\gamma):=\exp_H\circ\,\gamma$
is injective on a suitable
$0$-neighbourhood~$W$ in $\Hol_b(U_n,L(H))$
and a homomorphism of local groups
with respect to the Baker-Campbell-Hausdorff
multiplication on~$W$,
we deduce that the subgroup
$G_n$ of $\Hol(U_n,H)$
generated by $\Exp_n(\Hol_b(U_n,L(H)))$
can be made a Banach-Lie group
with Lie algebra
$\Hol_b(U_n,L(H))$.
The restriction map
$G_m\to G_n$, $\gamma\mto\gamma|_{U_n}$
is an injective, smooth homomorphism
for $n>m$,
and its differential $L(i_{n,m})\colon \Hol_b(U_m,L(H))\to\Hol_b(U_n,L(H))$,
$\gamma\mto\gamma|_{U_n}$ a compact
operator. Also, $G$
has a direct limit chart
(see \cite{GEM} and \cite{COM} for details).\\[4mm]
We remark that, for a more restrictive
class of Lie groups, there is a simple
criterion for the non-existence of small subgroups
(cf.\ \cite[Lemma~2.23]{GCX}
and \cite[Problem~II.5]{Nee}):\\[4mm]
{\bf Proposition.}
\emph{If a Lie group $G$
has an exponential map
which is a local homeomorphism at~$0$,
then $G$ does not have small torsion subgroups.
Also, $G$ does not have small subgroups
if} (\emph{and only if}) \emph{$L(G)$
admits a continuous norm.\Punkt}\\[3mm]
Combining Theorem~B\,(i)
and the preceding proposition,
we see that every Silva space $E=\bigcup_{n\in \N}E_n$
does not have small additive subgroups and hence
admits a continuous norm.
Since $\Gamma(K,H)$
has an exponential function which
is a local homeomorphism at~$0$ (see~\cite{GEM})
and $\Gamma(K,L(H))$ is a Silva space,
applying the proposition again
we get an alternative proof
for the non-existence of small
subgroups in $\Gamma(K,H)$.\\[2.5mm]
The preceding proposition
does not subsume Theorem~A
(although its hypotheses are satisfied by special cases
of direct limit groups as in \cite{KaM} or \cite{NRW}).
In fact, the exponential map of a direct limit group need
not be injective on any $0$-neighbourhood~\cite[Example~5.5]{DIR}.
\section{Some preliminaries concerning
direct limits}\label{secprelim}
Background information concerning direct limits
of topological groups,
topological spaces and Lie groups
can be found in \cite{DIR}, \cite{FUN}--\cite{Hir}
and \cite{TSH}.
We recall:
If $X_1\sub X_2\sub\cdots$ is an
ascending sequence of topological spaces
such that the inclusion maps $X_n\to X_{n+1}$
are continuous, then the final topology on $X:=\bigcup_{n\in \N} X_n$
with respect to the inclusion maps $X_n\to X$
makes $X$ the direct limit $\dl\,X_n$\vspace{-.6mm}
in the category of topological spaces and
continuous maps.
Thus, $S\sub X$ is open (resp., closed)
if and only if $S\cap X_n$ is open (resp., closed)
in~$X_n$ for each $n\in \N$.
If each $X_n$ is a locally convex real
topological vector space here
and each inclusion map $X_n\to X_{n+1}$
is continuous linear, then the \emph{locally
convex direct limit topology} on $X=\bigcup_{n\in \N}X_n$
is the finest locally convex vector topology making
each inclusion map $X_n\to X$ continuous
(see \cite{BTV}).
It is coarser then the direct limit topology,
and can be properly coarser.
For easy reference, let us compile
various well-known facts:
%
%
%\ma{baseDL}
\begin{la}\label{baseDL}
Let $X_1\sub X_2\sub\cdots$ be an ascending sequence
of topological spaces
and $X:=\bigcup_{n\in \N}X_n$,
equipped with the direct limit
topology.
\begin{itemize}
\item[\rm (a)]
If $S\sub X$ is open or closed,
then $X$ induces on~$S$
the topology making $S$ the direct limit
$S=\dl\,(S\cap X_n)$,\vspace{-.3mm}
where $S\cap X_n$ carries the topology induced
by~$X_n$.
\item[\rm (b)]
If $U_1\sub U_2\sub\cdots$
is an ascending sequence of open subsets $U_n\sub X_n$,
then $U:=\bigcup_{n\in \N}U_n$
is open in~$X$
and $U=\dl\,U_n$\vspace{-.7mm}
as a topological space.
\end{itemize}
\end{la}
\begin{proof}
(a) is immediate from the definition
of final topologies.
(b) is \cite[Lemma~1.7]{COM}.
\end{proof}
Given a topological space~$X$ and subset $Y\sub X$, we write
$Y^0$ for its interior.
A sequence $(U_k)_{k\in \N}$
of neighbourhoods of a point $x\in X$
is called a \emph{fundamental sequence}
if $U_k\supseteq U_{k+1}$ for each $k\in \N$
and $\{U_k\colon k\in \N\}$
is a basis of neighbourhoods for~$x$.
\section{\!\!Construction of neighbourhoods without subgroups}\label{secconstr}
The following lemma is the technical backbone
of our constructions.
In the lemma, $\cK$ denotes a set
of subsets of the given topological group~$G$,
with the following properties:
\begin{itemize}
\item[(a)]
$\cK$ is closed under finite unions; and
\item[(b)]
For each compact subset $K\sub G$,
the set
$\cK_K := \{S\in\cK\colon \mbox{$S$ is a neighbourhood
of $K$}\}$
is a basis of neighbourhoods of~$K$ in~$G$.
\end{itemize}
Of main interest are the three cases
where $\cK$ is, respectively,
the set of all closed subsets
of~$G$; the set of all compact subsets;
and the set of all subsets $S\sub G$ such that
$f(S)$ is compact, where
$f\colon G\to H$ is a given continuous homomorphism
to a topological group~$H$, such that each
$x\in G$ has a basis of neighbourhoods $U$
with compact image~$f(U)$.
%
%
%\ma{bbone}
\begin{la}\label{bbone}
Let $G$ be a topological group
without small subgroups
and $K\sub G$ be a compact set
that does not contain any
non-trivial subgroup of~$G$.
If $1\in K$,
then there exists a neighbourhood
$W$ of~$K$ in~$G$
which does not contain
any non-trivial subgroup of~$G$,
and such that $W\in \cK_K$.
Also, $W$ can be chosen
as a subset of any
given neighbourhood $X$ of~$K$.
\end{la}
\begin{proof}
We may assume that~$X$ is open.
Let $V\sub X$ be an open identity neighbourhood
such that $V$ does not contain any non-trivial
subgroup of~$G$, and $Q\sub V$ be a closed identity
neighbourhood of~$G$.
For each $x\in K\setminus Q^0$,
there exists $k\in \Z$ such that $x^k\not\in K$.
Let $J_x$ be a compact neighbourhood
of~$x$ in $K$ such that
$I_x:=\{y^k\colon y\in J_x\}\sub G\setminus K$.
Choose a closed neighbourhood
$P_x$ of~$I_x$ in $G\setminus K$
and let $A_x$ be a neighbourhood
of~$J_x$ in~$G$ such that $y^k\in P_x$
for each $y\in A_x$.
The set
$K\setminus Q^0$
being compact, we find
subsets $A_1,\ldots, A_m$ of~$G$
and compact subsets $J_1,\ldots, J_m$
of~$K$ such that
$K\setminus Q^0\sub\bigcup_{j=1}^m J_j$,
closed subsets $P_1,\ldots, P_m$
of~$G$ disjoint from~$K$
and $k_1,\ldots, k_m\in \Z$
such that
$J_j\sub A_j^0$ for each $j\in \{1,\ldots,m\}$
and $y^{k_j}\in P_j$ for each $y\in A_j$.
Then $P:=\bigcup_{j=1}^mP_j$
is a closed subset of~$G$
such that $P\cap K=\emptyset$.
After replacing $A_j$ with
a neighbourhood $\tilde{A}_j\in \cK_{J_j}$
of~$J_j$ contained in
$X\cap (A_j^0\setminus P)$
(which is an open neighbourhood of~$J_j$)
for each~$j$,
we may assume that $A:=\bigcup_{j=1}^m A_j$
and~$P$ are disjoint,
$A\sub X$, and $A\in \cK$.
Then $V\setminus P$ is a neighbourhood
of the compact set $Q\cap K$,
whence $B\sub V \setminus P$
for some $B\in \cK_{Q\cap K}$.
Then $W:=A\cup B\in \cK$.
We now show that $W$ does not contain
any non-trivial subgroup of~$G$.
Let $1\not=x\in W$.
Case~1: If $x\in A$, then $x^k\in P$
for some $k\in \Z$ and thus $x^k\not\in W$,
since $W$ and~$P$ are disjoint by construction.
Hence $\langle x\rangle\not\sub W$.
Case~2:
If $x\in B\sub V$, then $\langle x\rangle\not\sub V$,
whence there is $k\in \Z$ such that
$x^k\not\in V$.
If $x^k\in A$, then $\langle x^k\rangle\not\sub W$
by Case~1 and hence
$\langle x\rangle\not\sub W$ a fortiori.
If $x^k\not\in A$, then $x^k\not\in W$
(as $x^k\not\in B\sub V$ either)
and thus $\langle x\rangle \not\sub W$.
This completes the proof.
\end{proof}
%
%
%\ma{givecontrol}
\begin{rem}\label{givecontrol}
The proof of Lemma~\ref{bbone}
can easily be adapted
to get further information.
Namely, let
$C_1,\ldots, C_M$ be compact subsets
of~$K\setminus \{1\}$
and $\ell_1,\ldots,\ell_M$ be integers
such that $x^{\ell_j}\not\in K$
for each $j\in\{1,\ldots, M\}$
and $x\in C_j$.
Furthermore,
let $R,T\sub G$ be closed subsets
such that $T\cap K=\emptyset$ and $1\not\in R$.
We then easily achieve
that the following additional
requirements are met in the proof of Lemma~\ref{bbone}
(which will become vital later):
\begin{itemize}
\item[(a)]
$M\leq m$,
$C_j\sub A_j^0$ and $k_j=\ell_j$ for each $j\in \{1,\ldots, M\}$;
\item[(b)]
$W\cap T=\emptyset$ and $V\cap R=\emptyset$.
\end{itemize}
In fact,
we can simply
replace~$X$ by its intersection with
the open set $G\setminus T$
and choose~$V$ as a subset of $G\setminus R$
to ensure\,(b).
In the construction of
$k$, $J_x$, $P_x$ and $A_x$
described at the beginning of the proof of Lemma~\ref{bbone},
we can replace
$J_x$ with a compact neighbourhood~$J_j$
of $C_j$ in $K$ for $j\in \{1,\ldots, M\}$,
such that $I_j:=\{y^{\ell_j}\colon y\in J_j\}\sub G\setminus K$.
After enlarging the chosen finite cover of
$K\setminus Q^0$
by the preceding sets
if necessary,
we may assume that $m\geq M$
and $k_j=\ell_j$ as well as $C_j\sub A_j^0$,
for all $j\in \{1,\ldots, M\}$.
\end{rem}
\begin{rem}
It is a natural idea to try to prove, say,
Theorem~A for $G=\bigcup_{n\in \N}G_n$
in the following way:
Start with a compact identity neighbourhood
$W_1\sub G_1$ without non-trivial
subgroups,
and use Lemma~\ref{bbone}
recursively to obtain a sequence
$(W_n)_{n\in \N}$
of compact subsets $W_n\sub G_n$
such that $W_n$ has $W_{n-1}$
in its interior and does not contain
any non-trivial subgroup.
Then $W:=\bigcup_{n\in \N}W_n$ is
an identity neighbourhood in~$G$
and is a candidate
for an identity neighbourhood
not containing non-trivial subgroups.
But, unfortunately, it can happen
that~$W$ does contain non-trivial subgroups,
as the example $G_n:=\R^n$,
$G:=\R^{(\N)}=\dl\,G_n$,\vspace{-.9mm}
$W_n:=[{-n},n]^n$,
$W=\R^{(\N)}=G$ shows.
Therefore, this basic idea has to be refined,
and each $W_n$ has to be chosen in a
much more restrictive way.
The considerations from Remark~\ref{givecontrol}
will provide the required additional
control on the sets~$W_n$.
Further modifications
will be necessary to adapt the basic idea
to the (possibly) non-locally compact groups~$G_n$
in Theorem~B.
\end{rem}
\section{Proof of Theorem~B}
We start with several lemmas
which will help us to prove
Theorem~B.
The first lemma is a well-known fact from the theory of
Silva spaces,
but it is useful to recall its
proof here because details thereof
are essential for subsequent arguments.
%
%
%\ma{eqseq1}
\begin{la}\label{eqseq1}
Let $E_1\sub E_2\sub\cdots$
be an ascending sequence
of Banach spaces,
such that the inclusion map
$i_{n,m}\colon E_m\to E_n$ is a compact linear operator
whenever $n>m$.
Then there is an ascending sequence
$E_1\sub F_1\sub E_2\sub F_2\sub\cdots$
of Banach spaces with
continuous linear inclusion maps,
such that,
for each $n\in \N$,
there exists a norm $p_n$ on $F_n$
which defines the topology of~$F_n$
and has the property
that all closed $p_n$-balls $\wb{B}^{\,p_n}_r\!(x)$,
\emph{($r>0$, $x\in F_n$),}
are compact in $F_{n+1}$.
\end{la}
\begin{proof}
Let $B_n$ be the closed unit ball
in~$E_n$ with respect to some
norm defining its topology and
$K_n$ be the closure of $B_n$
in~$E_{n+1}$, which is compact
by hypothesis.
Let $F_n:=(E_{n+1})_{K_n}$
be the vector subspace of~$E_{n+1}$
spanned by~$K_n$
and~$p_n$ be the Minkowski
functional of~$K_n$ on~$F_n$.
Then $F_n$ is a Banach space,
by the corollary to Proposition~8
in \cite[Chapter~III, \S1, no.\,5]{BTV}.
The inclusion map
$F_n\to E_{n+1}$ is continuous,
and also the inclusion map
$E_n\to F_n$, since
$B_n\sub K_n=\wb{B}^{\,p_n}_1(0)$.
As $K_n$ is compact in~$E_{n+1}$
and the inclusion map
$E_{n+1}\to F_{n+1}$ is continuous,
$K_n$ is compact in $F_{n+1}$
(and hence also the image of any
ball $\wb{B}_r^{\,p_n}(x)$).
\end{proof}
%
%
%\ma{eqseq2}
\begin{la}\label{eqseq2}
If each $E_n$ is a Banach-Lie algebra
in the situation of Lemma~{\rm\ref{eqseq1}}
and each $i_{n,m}$ also is a
Lie algebra homomorphism,
then $F_n$ can be chosen as a Lie subalgebra
of~$E_{n+1}$ and it can be achieved
that $p_n$ makes $F_n$ a Banach-Lie algebra.
\end{la}
\begin{proof}
Since $[B_n,B_n]\sub r B_n$
for some $r>0$,
we have $[K_n,K_n]\sub r K_n$,
entailing that $F_n=\text{span}(K_n)$
is a Lie subalgebra of $E_{n+1}$
and the Lie bracket $F_n\times F_n\to F_n$
is a continuous bilinear map.
\end{proof}
Given a Banach-Lie group~$G$,
we let $\Ad^G\colon G\to\text{Aut}(L(G))$,
$x\mto\Ad_x^G$ be the adjoint homomorphism,
$\Ad^G_x:=L(c_x)$ with
$c_x\colon G\to G$, $c_x(y):=xyx^{-1}$.
%
%
%\ma{eqseq3}
\begin{la}\label{eqseq3}
Let $G_1\sub G_2\sub\cdots$
be an ascending sequence
of Banach-Lie groups,
such that the inclusion maps
$i_{n,m}\colon G_m\to G_n$
are smooth homomorphisms for
$n\geq m$ and $L(i_{n,m})\colon L(G_m)\to L(G_n)$
is a compact linear operator
whenever $n>m$.
Then there is an ascending sequence
$G_1\sub H_1\sub G_2\sub H_2\sub\cdots$
of Banach-Lie groups
such that, for each $n\in \N$,
there is a norm $p_n$ on $L(H_n)$
which defines the topology of~$L(H_n)$
and has the property that all closed $p_n$-balls $\wb{B}^{\,p_n}_r(x)$,
\emph{($r>0$, $x\in L(H_n)$),}
are compact in $L(H_{n+1})$.
\end{la}
\begin{proof}
We identify $L(G_n)$ with a Lie subalgebra
of~$L(G_{n+1})$ for each $n\in \N$.
By Lemma~\ref{eqseq2}, there is an ascending
sequence
\[
L(G_1)\; \sub \; F_1\; \sub \; L(G_2)\; \sub \; F_2\; \sub \;\cdots
\]
of Banach-Lie algebras
such that the inclusion maps are
continuous Lie algebra homomorphisms,
and such that, for each $n\in \N$,
there exists a norm $p_n$ on~$F_n$
defining its topology and
such that all closed $p_n$-balls
in~$F_n$ are compact subsets of $F_{n+1}$.
As in the proofs of Lemmas~\ref{eqseq1}
and~\ref{eqseq2}, we may assume
that the closed unit ball $K_n:=\wb{B}_1^{\,p_n}(0)$
of~$F_n$ is the closure in
$L(G_{n+1})$ of the closed unit ball~$B_n$
of $L(G_n)$.
We give
$S_n:=\langle \exp_{G_{n+1}}(F_n)\rangle$
the Banach-Lie group structure
making it an analytic subgroup of~$G_{n+1}$,
with Lie algebra~$F_n$.
For each $x\in G_n$, we have
\[
\Ad_x^{G_{n+1}}(B_n)\; =\; \Ad_x^{G_n}(B_n)\; \sub \; r B_n
\]
for some $r>0$,
whence $\Ad_x^{G_{n+1}}(K_n)\sub r K_n$
and hence $\Ad_x^{G_{n+1}}(F_n)\sub F_n$.
Note that the linear automorphism of $F_n$
induced by $\Ad^{G_{n+1}}_x$ is continuous,
by the penultimate inclusion.
As a consequence, the subgroup
$H_n:=\langle G_n\cup \exp_{G_{n+1}}(F_n)\rangle$
of $G_{n+1}$ can be given a
Banach-Lie group structure
with~$S_n$ as an open subgroup
(cf.\ Proposition~18 in
\cite[Chapter~III, \S1.9]{Bou}).
By construction, the
Banach-Lie groups $H_n$
have the desired properties.
\end{proof}
%
%
%\ma{nicebase}
\begin{la}\label{nicebase}
Let $f\colon G\to H$ be
a smooth homomorphism between
Banach-Lie groups such that,
for some norm~$p$ on $L(G)$ defining
its topology,
$L(f)\colon L(G)\to L(H)$
takes closed balls
in $L(G)$ to compact subsets
of~$L(H)$.
Then each $x\in G$ has a basis of closed
neighbourhoods~$U$ such that $f(U)$
is compact in~$H$.
Furthermore, every neighbourhood
of a compact subset $K\sub G$ contains
a closed neighbourhood~$A$
such that $f(A)$ is compact.
\end{la}
\begin{proof}
Since $G$ is a regular topological space
and $\exp_G$ a local homeomorphism at~$0$,
there is $R>0$ such that
$\exp_G|_{\wb{B}_R}$
is a homeomorphism onto its image
and $V_r:=\exp_G(\wb{B}_r)$
is closed in~$G$
for each $r\in \;]0,R]$,
where $\wb{B}_r:=\{x\in L(G)\colon p(x)\leq r\}$.
Exploiting the naturality of~$\exp$
and the hypothesis that $L(f).\wb{B}_r$
is compact in~$L(H)$,
we deduce that
$f(V_r)=f(\exp_G(\wb{B}_r))=\exp_H(L(f).\wb{B}_r)$
is compact in~$H$, for each $r\in \;]0,R]$.
Thus $\{V_r\colon r\in\;]0,R]\}$ is a basis
of closed neighbourhoods of~$1$ in~$G$
with compact image under~$f$.
Then
$\{x V_r\colon r\in\;]0,R]\}$
is a basis of closed neighbourhoods of $x\in G$
with compact image. The final assertion is an immediate
consequence.
\end{proof}
{\bf Proof of Theorem~B, assuming condition\,(i).}
We define Banach-Lie groups $H_n$
as in Lemma~\ref{eqseq3}.
After replacing $G_n$ with~$H_n$
for each $n\in \N$,
we may assume that
each point in~$G_n$ has a basis
of neighbourhoods in~$G_n$
which are compact in~$G_{n+1}$
(see Lemma~\ref{nicebase}).
We now construct, for each $n\in \N$:
\begin{itemize}
\item An
identity neighbourhood
$W_n\sub G_n$
such that $W_n$,
when considered as a subset~$K_n$
of $G_{n+1}$, becomes compact;
\item
A fundamental sequence $(Y_k^{(n)})_{k\in \N}$
of open identity neighbourhoods
in~$K_n$;
\item
For some $m_n\in \N_0$,
a family $(C_j^{(n)})_{j=1}^{m_n}$
of subsets $C_j^{(n)}$ of $W_n\setminus \{1\}$
which are compact in~$G_{n+1}$; and
\item
A function $\kappa_n\colon \{1,\ldots, m_n\} \to \Z$,
\end{itemize}
with the following properties:
\begin{itemize}
\item[(a)]
If $n>1$, then
$W_{n-1}$ is contained in the interior $W_n^0$
of $W_n$ relative~$G_n$;
\item[(b)]
$W_n$ does not contain
any non-trivial subgroup of~$G_n$;
\item[(c)]
For each $j\in \{1,\ldots, m_n\}$ and $x\in C_j^{(n)}$,
we have $x^{\kappa_n(j)}\not\in W_n$;
\item[(d)]
If $n>1$, then $m_n\geq m_{n-1}$
and $C_j^{(n-1)}\sub C_j^{(n)}$
as well as $\kappa_n(j)=\kappa_{n-1}(j)$,
for all $j\in \{1,\ldots, m_{n-1}\}$;
\item[(e)]
For all positive integers $\ell<n$,
we have $K_\ell\setminus Y_n^{(\ell)}\sub
\bigcup_{j=1}^{m_n} C_j^{(n)}$.
\end{itemize}
If this construction is possible,
then $U:=\bigcup_{n\in \N}W_n=\bigcup_{n\geq 2}W_n^0$
is an open identity
neighbourhood in $G=\dl\,G_n$,\vspace{-.3mm}
using\,(a) and Lemma~\ref{baseDL}\,(b).
Furthermore, $U$ does not contain any
non-trivial subgroup of~$G$.
In fact: If $1\not=x\in U$,
there is $m\in \N$ such that $x\in W_m$.
Then $x\in K_m\setminus Y^{(m)}_n$
for some $n>m$,
and thus $x\in C_j^{(n)}$ for some $j\in \{1,\ldots, m_n\}$,
by\,(e). By (c) and (d),
we have $x^{\kappa_n(j)}\not\in W_k$
for each $k\geq n$,
whence $x^{\kappa_n(j)}\not\in U$
and thus $\langle x\rangle \not\sub U$.\\[2.8mm]
It remains to carry out the construction.
Since $G_1$ is a Banach-Lie group,
it does not have small subgroups,
whence we find an identity neighbourhood
$W_1$ in~$G_1$ which does not contain
any non-trivial subgroup of~$G_1$.
By Lemma~\ref{nicebase},
after replacing~$W_1$ be a smaller
identity neighbourhood, we may assume that
$W_1$, considered as subset~$K_1$ of $G_2$,
becomes compact.
We set $m_1:=0$,
$\kappa_1:=\emptyset$,
and choose any fundamental sequence
$(Y^{(1)}_k)_{k\in \N}$
of open identity neighbourhoods
of~$K_1$, which is possible because~$G_2$
(and hence~$K_1$) is metrizable.\\[2.8mm]
Let $N$ be an integer $\geq 2$ now
and suppose that $W_n$,
$(Y^{(n)}_k)_{k\in \N}$,
$(C_j^{(n)})_{j=1}^{m_n}$
and $\kappa_n$  have been constructed
for $n\in\{1,\ldots, N-1\}$,
such that (a)--(e) hold.
Then $R:=\bigcup_{\ell<N}K_\ell\setminus Y_N^{(\ell)}$
and $T:=\bigcup_{j=1}^{m_{N-1}}\{x^{\kappa_{N-1}(j)}\colon x\in C_j^{(N-1)}\}$
are compact subsets of $G_N$
such that $1\not\in R$
and $T\cap W_{N-1}=\emptyset$.
We now apply Lemma~\ref{bbone}
to $G_N$ and its
compact subset $K:=K_{N-1}$,
with $\cK$ the set of all subsets of~$G_N$
which are compact in $G_{N+1}$.
Let $A_1,\ldots,A_m$,
$k_1,\ldots, k_m$,
$V$, $A$ and $W_N:=W\in \cK$ be as described
in Lemma~\ref{bbone} and its proof.
As explained in Remark~\ref{givecontrol},
we may assume that $V\cap R =\emptyset$,
$W\cap T=\emptyset$,
$m\geq m_{N-1}$, $C_j^{(N-1)}\sub A_j^0$
for each $j\in \{1,\ldots, m_{N-1}\}$,
and $k_j=\kappa_{N-1}(j)$.
Set $m_N:=m$,  $C_j^{(N)}:=A_j$
for $j\in \{1,\ldots, m_N\}$,
and $\kappa_N(j):=k_j$.
Let $(Y^{(N)}_k)_{k\in\N}$
be any fundamental sequence
of open identity neighbourhoods
in $K_N:=W_N$, considered
as a compact subset of~$G_{N+1}$.
If $\ell<N$, then $K_\ell\setminus Y_N^{(\ell)}\sub R$
and hence $(K_\ell\setminus Y_N^{(\ell)})\cap V=\emptyset$,
entailing that $K_\ell\setminus Y_N^{(\ell)}\sub
W\setminus V\sub A=\bigcup_{j=1}^{m_N}C_j^{(N)}$.
Thus (a)--(e) hold
for all $n\in \{1,\ldots, N\}$.\pagebreak

\noindent
{\bf Proof of Theorem~B, assuming condition\,(ii).}
Let $\phi\colon \tilde{Z} \to \tilde{H}\sub L(G)$
be a direct limit chart of~$G$ around~$1$,
such that $\phi(1)=0$.
Thus $\tilde{Z}=\bigcup_{n\in \N}Z_n$,
$\tilde{H}=\bigcup_{n\in \N}H_n$,
and $\phi=\bigcup_{n\in \N}\phi_n$
for certain charts $\phi_n\colon Z_n\to H_n$
of~$G_n$, such that $Z_n\sub Z_{n+1}$
and $\phi_{n+1}|_{Z_n}=\phi_n$
for each $n\in \N$.
By \cite[Proposition~7.12]{GaG},
$L(G)$ is a $k_\omega$-space
and $L(G)=\dl\,L(G_n)$\vspace{-.5mm}
also as a topological space.
By \cite[Proposition~4.2\,(g)]{GaG},
$H_1$ has an open $0$-neighbourhood
$V_1$ which is a $k_\omega$-space.
By Proposition~4.2\,(g) and Lemma~4.3 in \cite{GaG},
$V_1$ has an open neighbourhood
$V_2$ in $H_2$ which is
a $k_\omega$-space.
Proceeding in this way, we find
an ascending sequence $V_1\sub V_2\sub\cdots$
of open $0$-neighbourhoods
$V_n\sub H_n$,
such that each $V_n$ is a $k_\omega$-space.
By Lemma~\ref{baseDL}\,(b),
$V:=\bigcup_{n\in \N}V_n\sub H$
is open in $L(G)$
and $V=\dl\, V_n$\vspace{-.5mm}
as a topological space,
whence $V$ is a $k_\omega$-space
by \cite[Proposition~4.5]{GaG}.
For each $j\in \N$,
choose a $k_\omega$-sequence
$(L_n^{(j)})_{n\in \N}$
for $V_j$.
We may assume that $0\in L_1^{(1)}$.
After replacing $L_n^{(j)}$ with
$\bigcup_{i=1}^j L_n^{(i)}$,
we may assume that $L_n^{(i)}\sub L_n^{(j)}$
for all positive integers $i\leq j$ and~$n$.
Then $L_n^{(n)}$ is a $k_\omega$-sequence
for~$V$ (see the first half of the proof
of Proposition~4.5 in~\cite{GaG}),
and thus $K_n:=\phi^{-1}(L_n^{(n)})$
defines a $k_\omega$-sequence $(K_n)_{n\in \N}$ for
the open identity neighbourhood
$Z:=\phi^{-1}(V)\sub G$.
Note that $K_n=\phi_n^{-1}(L_n^{(n)})$
is a compact subset of~$G_n$,
and $1\in K_1$.
Because $L(G_n)$ admits a continuous
norm, the compact set $L_n^{(n)}$
is metrizable and hence also~$K_n$.
We now construct, for each $n\in \N$:
\begin{itemize}
\item
A compact identity neighbourhood $W_n$ in $K_n$;
\item
A fundamental sequence $(Y_k^{(n)})_{k\in \N}$
of open identity neighbourhoods
in~$W_n$;
\item
For some $m_n\in \N_0$,
a family $(C^{(n)}_j)_{j=1}^{m_n}$
of subsets $C_j^{(n)}\sub W_n\setminus \{1\}$,
and a function $\kappa_n\colon \{1,\ldots, m_n\} \to \Z$,
\end{itemize}
such that conditions (b)--(e) from the proof of Theorem~B\,(i)
are satisfied and also
\begin{itemize}
\item[(a)$'$]
If $n>1$, then
$W_{n-1}$ is contained in the interior $W_n^0$ of $W_n$
relative~$K_n$.
\end{itemize}
If this construction is possible,
then $U:=\bigcup_{n\in \N}W_n=\bigcup_{n\geq 2}W_n^0$
is an open identity
neighbourhood in $Z=\dl\,K_n$\vspace{-.5mm}
(by (a)$'$ and Lemma~\ref{baseDL}\,(b)),
and hence in~$G$.
Furthermore, $U$ does not contain any
non-trivial subgroup of~$G$,
by the same argument as above.\\[3mm]
To carry out the construction,
we recall first that as $G_n$ has an exponential map
which is a local homeomorphism at~$0$
and $L(G_n)$ admits a continuous
norm, $G_n$ does not have small subgroups
(by the proposition in the Introduction).
In particular,
we find a closed identity neighbourhood
$\tilde{W}_1$ in~$G_1$ which does not contain
any non-trivial subgroup of~$G_1$.
Then $W_1:=\tilde{W}_1\cap K_1$
is a compact identity neighbourhood
in~$K_1$.
We set $m_1:=0$,
$\kappa_1:=\emptyset$,
and choose any fundamental sequence
$(Y^{(1)}_k)_{k\in \N}$
of open identity neighbourhoods
of~$W_1$ (which is possible
because~$K_1$ is metrizable).\\[2.8mm]
Let $N$ be an integer $\geq 2$ now
and suppose that $W_n$,
$(Y^{(n)}_k)_{k\in \N}$,
$(C_j^{(n)})_{j=1}^{m_n}$
and $\kappa_n$ have been constructed
for $n\in\{1,\ldots, N-1\}$ such that (a)$'$
and (b)--(e) hold.
Then $R:=\bigcup_{\ell<N}K_\ell\setminus Y_N^{(\ell)}$
and $T:=\bigcup_{j=1}^{m_{N-1}}\{x^{\kappa_{N-1}(j)}\colon x\in C_j^{(N-1)}\}$
are compact subsets of $G_N$
such that $1\not\in R$
and $T\cap W_{N-1}=\emptyset$.
We now apply Lemma~\ref{bbone}
to $G_N$ and its
compact subset $K:=K_{N-1}$,
with $\cK$ the set of all closed
subsets of~$G_N$.
Let $A_1,\ldots,A_m$,
$k_1,\ldots, k_m$,
$V$, $A$ and $W\in \cK$ be as described
in Lemma~\ref{bbone} and its proof.
As explained in Remark~\ref{givecontrol},
we may assume that $V\cap R =\emptyset$,
$W\cap T=\emptyset$,
$m\geq m_{N-1}$, $C_j^{(N-1)} \sub A_j^0$
for each $j\in \{1,\ldots, m_{N-1}\}$,
and $k_j=\kappa_{N-1}(j)$.
Set $W_N:=W\cap K_N$,
$m_N:=m$,
$C_j^{(N)}:=A_j\cap K_N$
for $j\in \{1,\ldots, m\}$,
and $\kappa_N(j):=k_j$.
Let $(Y^{(N)}_k)_{k\in\N}$
be any fundamental sequence
of open identity neighbourhoods
in~$W_N$.
If $\ell<N$, then $K_\ell\setminus Y_N^{(\ell)}\sub R$
and hence $(K_\ell\setminus Y_N^{(\ell)})\cap V=\emptyset$,
entailing that $K_\ell\setminus Y_N^{(\ell)}\sub
K_N\cap (W\setminus V)\sub K_N\cap A=\bigcup_{j=1}^{m_N}C_j^{(N)}$.
Thus (a)$'$ and (b)-(e)
hold for all $n\in \{1,\ldots, N\}$.\vspace{3mm}\Punkt

\noindent
\emph{Acknowledgement.}
The author thanks K.-H. Neeb (Darmstadt)
for an inspiring question.\vspace{-4mm}
\noindent
{\footnotesize{\bf Helge Gl\"{o}ckner},
\,TU~Darmstadt, FB~Mathematik~AG~5,
Schlossgartenstr.\,7,
64289 Darmstadt, Germany\\[.2mm]
E-Mail:
\,{\tt gloeckner@mathematik.tu-darmstadt.de}}

{\footnotesize
\begin{thebibliography}{WW}\itemsep=-.03pc
%
%
% 
\bibitem{Aus} Au\ss{}enhofer, L.,
{\em Contributions to the duality theory of Abelian
topological groups and to the theory
of nuclear groups}, Dissertationes Math.\ {\bf 384},
1999.
%
%
\bibitem{BTV}
Bourbaki, N., ``Topological Vector Spaces,
Chapters~1--5,'' Springer,
Berlin, 1987.
%
%
\bibitem{Bou}
Bourbaki, N., ``Lie Groups and Lie Algebras,
Chapters~1--3,''
Springer, Berlin, 1989.
%
%
\bibitem{Flo}
Floret, K., \emph{Lokalkonvexe Sequenzen mit kompakten
Abbildungen}, J. Reine Angew.\ Math.\ {\bf 247} (1971),
155--195.
%
%
\bibitem{GCX} Gl\"{o}ckner, H.,
{\em Lie group structures on quotient groups
and universal complexifications for infinite-dimensional
Lie groups}, J. Funct.\ Anal.\ {\bf 194} (2002), 347--409.
%
%
\bibitem{DIR} Gl\"{o}ckner, H.
\emph{Direct limit Lie groups and manifolds},
J. Math.\ Kyoto Univ.\
{\bf 43} (2003), 1--26.
%
%
\bibitem{GEM} Gl\"{o}ckner, H.,
\emph{Lie groups of germs of analytic mappings},
pp.\ 1--16 in: T. Wurzbacher (Ed.),
``Infinite Dimensional Groups and Manifolds,''
IRMA Lecture Notes in Math.\ and Theor.\ Physics,
de Gruyter, 2004.
%
%
\bibitem{FUN} Gl\"{o}ckner, H.
\emph{Fundamentals of direct limit Lie theory},
Compos.\ Math.\ {\bf 141} (2005),
1551--1577.
%
%
\bibitem{COM} Gl\"{o}ckner, H.,
\emph{Direct limits of infinite-dimensional
Lie groups compared to direct limits
in related categories}, arXiv:math.GR/0606078.
%
%
\bibitem{GaG} Gl\"{o}ckner, H. and
R. Gramlich,
\emph{Final group topologies, Phan systems
and Pontryagin duality},
preprint, arXiv:math.GR/0603537.
%
%
%
\bibitem{Han}
Hansen, V.\,L.,
{\em Some theorems on direct limits
of expanding systems of manifolds},
Math.\ Scand.\ \textbf{29} (1971), 5--36.
%
%
\bibitem{Hir}
Hirai, T., H. Shimomura,
N. Tatsuuma and E. Hirai,
{\em Inductive limits
of topologies, their direct products, and problems
related to algebraic structures},
J. Math.\ Kyoto Univ.\ \textbf{41} (2001),
475--505.
%
%
\bibitem{KaM}
Kriegl, A. and P.\,W. Michor,
``The Convenient Setting of
Global Ana\-lysis,''
AMS, 1997.
%
%
\bibitem{Mil}
Milnor, J.,
{\it Remarks on infinite-dimensional Lie groups},
pp.\ 1007--1057 in:
``Relativit\'{e}, Groupes et Topologie II,''
B. DeWitt and R. Stora (Eds),
North-Holland, Amsterdam,
1984.
%
%
\bibitem{NRW}
Natarajan, L., E. Rodr\'{\i}guez-Carrington
and J.\,A. Wolf,
{\it Differentiable structure for direct limit groups},
Letters Math.\ Phys.\ {\bf 23} (1991), 99--109.
%
%
\bibitem{Nee}
Neeb, K.-H.,
\emph{Towards a Lie theory of locally convex groups},
TU Darmstadt Preprint {\bf 2459}, 2006.
%
%
\bibitem{TSH}
Tatsuuma, N., H. Shimomura, and T. Hirai,
\emph{On group topologies and unitary representations of inductive
limits of topological groups and the case of
the group of diffeomorphisms},
J. Math.\ Kyoto Univ.\ {\bf 38} (1998), 551--578.\vspace{1mm}
%
%
\end{thebibliography}}
\end{document}